\def\AJM{Amer.\ J.\ Math.}

\def\AP{Ann. Probab.}
\def\AU{Auburn University}

\def\emph{\em}

\def\LNM{Lecture Notes in Math.}

\def\PMS{Probab.\ Math.\ Statist.}
\def\PTRF{Probab.\ Theory Related Fields}

\def\SM{Studia Math.}
\def\Spr{Springer Vg.}

\def\ll{\Big\langle\, }
\def\rr{\,\Big\rangle}

\def\Sum{ \displaystyle\sum }
\def\Frac{ \displaystyle\frac }

\def\Rec#1{ \Frac{1}{#1} }

\def\df{\buildrel{\sf df}\over=}
\def\relbuild#1\under#2{\mathrel{\mathop{#2}\limits_{#1}}}

\def\D{\buildrel {\cal D}\over =}
\def\frac#1#2{{#1\over#2}}

\def\rec#1{ \frac{1}{#1} }

\def\I#1{1\hspace{-3pt} {\rm I}_{#1}}
\def\bm#1{\mbox{\boldmath $#1$}}
\def\ca#1{{\cal #1}}
\def\ov#1{\overline{#1}}

\def\wt#1{\widetilde{#1}}

\def\QED{\hfill \rule{2.5mm}{2.5mm}
\vspace{7pt}

}

\def\L{{\Bbb L}}
\def\BE{{\Bbb E}}

\def\R{{\Bbb R}}
\def\N{{\Bbb N}}
\def\E{\mbox{\sf E}}

\def\P{\mbox{\sf P}}
\def\eps{\varepsilon}
\def\beps{\bm\eps}

\def\refrm#1{{\rm(\ref{#1})}}
\def\set#1{\left\{\,#1\,\right\} }
\def\<{\langle}
\def\>{\rangle}

\def\Proof{{\sc Proof.~~}}

\def\otk{{\otimes k}}

\def\BX{{\Bbb X}}
\def\wh#1{\widehat{#1}}

\def\AJM{Amer.\ J.\ Math.}

\def\AP{Ann. Probab.}
\def\AU{Auburn University}

\def\LNM{Lecture Notes in Math.}

\def\PMS{Probab.\ Math.\ Statist.}
\def\PTRF{Probab.\ Theory Related Fields}
\def\SM{Studia Math.}

\def\Spr{Springer Vg.}

\def\bi{\bm i}

\documentstyle[11pt,leqno]{article}

\expandafter\chardef\csname pre amssym.def at\endcsname=\the\catcode`\@
\catcode`\@=11

\def\undefine#1{\let#1\undefined}
\def\newsymbol#1#2#3#4#5{\let\next@\relax
 \ifnum#2=\@ne\let\next@\msafam@\else
 \ifnum#2=\tw@\let\next@\msbfam@\fi\fi
 \mathchardef#1="#3\next@#4#5}
\def\mathhexbox@#1#2#3{\relax
 \ifmmode\mathpalette{}{\m@th\mathchar"#1#2#3}%
 \else\leavevmode\hbox{$\m@th\mathchar"#1#2#3$}\fi}
\def\hexnumber@#1{\ifcase#1 0\or 1\or 2\or 3\or 4\or 5\or 6\or 7\or 8\or
 9\or A\or B\or C\or D\or E\or F\fi}

\font\tenmsa=msam10
\font\sevenmsa=msam7
\font\fivemsa=msam5
\newfam\msafam
\textfont\msafam=\tenmsa
\scriptfont\msafam=\sevenmsa
\scriptscriptfont\msafam=\fivemsa
\edef\msafam@{\hexnumber@\msafam}
\mathchardef\dabar@"0\msafam@39

\font\tenmsb=msbm10
\font\sevenmsb=msbm7
\font\fivemsb=msbm5
\newfam\msbfam
\textfont\msbfam=\tenmsb
\scriptfont\msbfam=\sevenmsb
\scriptscriptfont\msbfam=\fivemsb
\edef\msbfam@{\hexnumber@\msbfam}

\def\Bbb#1{\fam\msbfam\relax#1}
\font\teneufm=eufm10
\font\seveneufm=eufm7
\font\fiveeufm=eufm5
\newfam\eufmfam
\textfont\eufmfam=\teneufm
\scriptfont\eufmfam=\seveneufm
\scriptscriptfont\eufmfam=\fiveeufm

\catcode`\@=\csname pre amssym.def at\endcsname

\newtheorem{theorem}{\sc Theorem}[section]
\newtheorem{lemma}[theorem]{\sc Lemma}
\newtheorem{proposition}[theorem]{\sc Proposition}
\newtheorem{remark}{\sc Remark}

\newcommand{\footer}[1]{{\def\thefootnote{}\footnotetext{#1}}}

\makeatletter
\def\EquationsBySection{\def\theequation{\thesection.\arabic{equation}}
\@addtoreset{equation}{section}}
\makeatother
\EquationsBySection
\renewcommand{\baselinestretch}{1.2}

\def\lnorm{\|}
\def\rnorm{\|}

\def\otk{{\otimes k}}

\def\BX{{\Bbb X}}
\def\L{{\Bbb L}}
\def\wh#1{\widehat{#1}}

\def\wt#1{\widetilde{#1}}

\def\dl{{\delta}}
\def\ddl{{\bm\delta}}

\def\ov{\overline}

\def\l({\bm{\big(\hspace{-3pt} \big(}\, }
\def\r){\,\bm{\big)\hspace{-3pt} \big)} }
\def\ll{\Big\langle\, }
\def\rr{\,\Big\rangle}

\begin{document}

\footer{{\bf Key words and phrases}: decoupling principle, symmetric tensor
products, random polynomials, multiple random series, multiple stochastic
integrals, random multilinear forms, random chaos, tail
inequalities, polarization, symmetrization, Banach space,
rearrangement invariant, Orlicz space, Lorentz space, Rademacher sequence, 
U-statistics}
\footer{{\bf AMS (1980) Classification} {\em  Primary}: 60B11, 46M05 {\em
Secondary}: 60H07, 46E30, 60E15, 62H05,  62G30}
\vspace{10pt}

\begin{center}
{\large\bf CONTRACTION AND DECOUPLING INEQUALITIES FOR 
MULTILINEAR FORMS AND U-STATISTICS}
\vspace{5pt}

{\sc V.H.~de la Pe\~na}$^{\dag}$
\footer
{$^{\dag}$ Research supported in part by N.S.F.\ Grant DMS 91-08006 and a grant
from the CNRS of France.}

{\small\it Department of Statistics, Columbia University,

New York, NY 10027.}
\vspace{5pt}

{\sc S.J.~Montgomery-Smith}$^{\ddag}$
\footer
{$^{\ddag}$ Research supported in part by N.S.F.\ Grants DMS 9001796 and DMS
9201357.}%

{\small\it 
Department of Mathematics, University of Missouri,

Columbia, MO 65211.}
\vspace{5pt}

{\sc Jerzy Szulga} 

{\small\it Department of Mathematics, Auburn University,

Auburn, AL 36849
}
\end{center}
\vspace{7pt}

\parbox{6in}{\small We prove decoupling inequalities for 
random polynomials in independent random variables with coefficients in vector 
space. We use various means of comparison, including rearrangement invariant
norms (e.g., Orlicz and Lorentz norms), tail distributions, tightness,
hypercontractivity, etc. }

\renewcommand{\baselinestretch}{0.9}
{\small \tableofcontents }
\eject

\renewcommand{\baselinestretch}{1.2}
\vspace{7pt}

\section{Introduction}
\subsection{Background and scope of the paper}

Decoupling principles stem from the theory of martingale transforms (cf.
\cite{Bur:mf}). For homogeneous random forms of rank $k\ge 2$, decoupling
principles were introduced in \cite{McCT:realdec,Kwa:dec,McCT:bandec} (in some
special cases, they were known to Pisier, cf. \cite{McCT:dyadic}), and
subsequently became essential tools in multiple integration (cf.
\cite{KalS,Szu:Umli,KwaW:double,RosW:clock,RST}.
One of the most appealing interpretations of such an principle is the
reducibility of a study of multiple random series (respectively, of multiple
stochastic integrals) to a consecutive treatment of single random series
(respectively, of single stochastic integrals which allows one to treat a
multiple integral as an It\^o-type iterate integral).  The concept of a random
chaos goes back to N.\ Wiener \cite{Wie:chaos} (see also \cite{WieW}), who
elaborated what we call here a real-valued coupled Gaussian chaos. Decoupling
inequalities may be viewed as embedding-projection procedures, since  a
decoupled random chaos is nothing but a lacunary random chaos.  In comparison
to the classical $L^2$-theory of multiple summation or integration, decoupling
principles make up the lack of $L^2$-isometries.  

Since the first
publication of the aforementioned decoupling principle, the theory has branched
into several directions. For example, comparison of tangent processes (cf.
\cite{Hit:tan,Pen:Wald,Pen:bound},) is akin to the classical decoupling 
principle.  Further
contributions can be found, e.g.,  in
\cite{Kwa:dec,DeA:dec,Hit:tan,Zin:comp,NolP,Pen}.  In some of the
aforementioned papers (e.g., \cite{McCT:realdec,McCT:bandec,Kwa:dec,DeA:dec})
the symmetry of random variables is essential for the fulfillment of the
decoupling principle, while other papers (e.g.,
\cite{Hit:tan,Zin:comp,KalS}) point out the role of positivity.
Norms of $L^p$-spaces, or more general, of Orlicz spaces (basically, subject to
growth restriction) have become typical means of comparison of two classes of
vector random variables. 

In this paper, we prove decoupling inequalities for random variables that are
not necessarily symmetric. Theorems~\ref{AB} and~\ref{F} in Section \ref{rob}, 
and Theorem~\ref{ABtail} in Section 3, are our main results.
The decoupling principle by means of probability tails, 
Theorem \ref{ABtail}, immediately
ensures the parity of tightness of two types of chaoses (that Gaussian
decoupled and coupled chaoses are simultaneously tight was proved in
\cite{Kwa:dec}). 

A number of decoupling results are obtained for arbitrary rearrangement
invariant norms and Orlicz functionals. 
In particular, we provide one extended example regarding certain
Lorentz norms (important in the approximation theory).
Another application is the decoupling principle for U-statistics (a
result as in Theorem \ref{F} was  proven in \cite{Pen}).

The utilized techniques are based on ideas, borrowed from
\cite{Kwa:dec}, while some are taken from \cite{KwaW:book}.
Proofs are straightforward and point out the algebraic nature of decoupling
that is fruitfully merged with a widely understood context of convexity.   A
rule of thumb is that, in the field of random diagonal-free polynomials, a
``definable'' is ``decouplable''.  The obtained robust
constants are tightly estimated, and are sharper than constants known 
before.

In the last section, we show tail probability decoupling results for
polynomials of symmetric random variables.  This section makes use of 
techniques from \cite{AsmM}.

\subsection{Notation}
Random variables in this paper are defined on
a separable probability space $(\Omega, \ca A,\P)$ that is rich enough to carry
independent sequences.  A sequence of real random variables is
denoted by $\bm\xi = (\xi_1,\xi_2,\dots)$, 
and a matrix of real random variables is denoted by
$\BX=[\bm\xi_1,\ldots,\bm\xi_n]$, where $\bm\xi_j = (\xi_{j1},\xi_{j2},\ldots)$.
We will make particular
use of one particular sequence, the Rademacher
sequence 
$\beps=(\eps_i)$, where $\eps_i$ are independent random variables taking
values $\pm 1$ with probability $1/2$.   

Let
$\BE=(\BE,\|\cdot\|)$ denote a real Banach space.
We will be considering $\BE$\ valued random variables, that is, 
strongly measurable mappings from $\Omega$
into $\BE$. 

\def\ik{{i_1,\dots,i_k}}
Let $k$ be a positive integer and $f = (f_\ik)$  be an array of vectors
from $\BE$\ taking only finitely many non-zero values.
Throughout the paper, all such arrays are assumed to
vanish on diagonals (We will say ``{\em diagonal-free}''), i.e., 
$f_\ik=0$, if
at least two indices $i_j,i_{j'}$ are equal. 

The main object of our interest will be the following
$k$-homogeneous random polynomial 
\[
Q(f;\BX) \df Q(f;\bm\xi_1,\ldots,\bm\xi_k)\df
\sum_{\ik} 
f_\ik \xi_{1 i_1} \cdots \xi_{k i_k}.
\]
We will be desiring to compare this random polynomial with the `undecoupled'
version, that is
\begin{equation}\label{undecoupled}
Q(f; \bm\xi^k) \df
Q(f;\bm\xi,\ldots,\bm\xi)\df
\sum_{\ik} 
f_\ik \xi_{i_1} \cdots \xi_{i_k}.
\end{equation}
The first term in the above definition will be introduced as a notational
convenience.  We will be quite free to stretch the use of this notation.  
So for
example, we might write
\[
Q(f; \bm\xi^r , \bm\eta^{k-r})
=
Q(f; \bm\xi,\ldots,\bm\xi , \bm\eta,\ldots,\bm\eta)
=
\sum_{\ik} 
f_\ik \xi_{i_1} \cdots \xi_{i_r} \eta_{i_{r+1}} \cdots \eta_{i_k}.
\]
We hope
to convince the reader of the value of this
notation, as it enables us to write many of the proofs in a more
compact form, and may ultimately lead to a clearer thinking on the
subject.  For the unconvinced reader, we hope that we have explained
the notation sufficiently that he will be able to rewrite all the
proofs and statements below in a more familiar form.

Many of the inequalities that we introduce require the array $f$\ to satisfy
certain symmetry conditions, and so for this reason we introduce the
symmetrized version of $f$:
\[\wh{f}_\ik
\df\rec{k!}\sum_\sigma f_{i_{\sigma_1},\ldots,i_{\sigma_k}}
\]
where the sum is taken over all permutations of the set
$[1,k]=\set{1,\ldots,k}$, and let $\wh{Q}(f;\cdot)=Q(\wh{f},\cdot)$. Note that
for the undecoupled random polynomial symmetry makes no change: 
$\wh{Q}(f;\bm\xi^k)={Q}(f;\bm\xi^k)$.

In the sequel, we will occasionally refer to tetrahedral arrays, i.e. 
$f$ such that $f_\ik = 0$, if indices fail to satisfy  $i_1<\ldots<i_k$.

We will make frequent use of
the following identity, which is known as the Mazur-Orlicz polarization 
formula \cite{MazO}:
\begin{equation}\label{MO}
\wh{Q}(f;\bm\xi_1,\ldots,\bm\xi_k)=
\rec{k!}\Sum_{\ddl=(\dl_1,\ldots,\dl_k)\in \set{0,1}^k}
(-1)^{k-|\ddl|}Q(f;(\dl_1\bm\xi_1+\ldots+\dl_k\bm\xi_k)^k),
\end{equation}
where $|\bm\delta|=\sum_i \delta_i$. Switching to a 
Rademacher sequence $\beps$, we can write
\begin{equation}\label{MORad}
\wh{Q}(f,\bm\xi_1,\ldots,\bm\xi_k)=
\Rec{k!}
\E\, 
\eps_1\cdots \eps_k \,Q(f,(\sum_{i=1}^k\eps_i \bm\xi_i)^{k}),
\end{equation}
where the expectation is only over the Rademacher sequence $\bm\eps$.

\paragraph{Rearrangement invariant spaces}
By $(\L,\|\cdot\|_\L)$  we denote a rearrangement invariant Banach space of
integrable random variables (so that the norm of a random variable depends
only on its probability distribution), $\L\subset L^1(\P)$, defined on a
separable probability space $(\Omega, \ca A, \P)$ that is rich enough to carry
independent sequences.
For more information on  rearrangement invariant spaces, we refer the reader
to, e.g., \cite{LinT}.
The basic examples of rearrangement invariant spaces are $\L = L_p$,
for $1 \le p \le \infty$, (nothing more is needed in many parts of this paper),
Orlicz spaces and Lorentz spaces.  We will sometimes use the abbreviation `r.i.'
for `rearrangement invariant.'

The important property of rearrangement invariant spaces that we shall use
is the following:
\begin{equation}\label{L}
\mbox{\em Conditional expectations are contractions acting on $\L$}.
\end{equation}
The reader unfamiliar with r.i.\ spaces should note that this is true of
$L_p$.

We denote by $\L(\BE)$ the Banach space of $\BE$-valued random variables (i.e.,
strongly measurable mappings from $\Omega$ into $\BE$) whose norms belong to
$\L$, and let $\|\theta\|_{\L(\BE)}=\|\,\|\theta\|_{\BE}\,\|_{\L}$.  
Thus if $\L = L_p$, then $\|\theta\|_{\L(\BE)} = \|\theta\|_{L_p(\BE)}
= (\E \|\theta\|_\BE ^p)^{1/p}$.
In
the sequel we sometimes omit the subscript indicating the space, if it 
causes no ambiguity.  

\section{Decoupling for r.i.\ norms}
\label{rob}
\paragraph{Interchangeability}
In the sequel, we will use several times the following elementary feature of
interchangeable random sequences $ \bm\xi_1,\ldots,\bm\xi_r$ (i.e., such that
each permutation has the same distribution). 
Suppose that each $\bm\xi_k$\ is itself a sequence of
independent random variables.
Denote by $\ca G_r$ the
$\sigma$-field spanned by  $\sum_{j=1}^r
\bm\xi_j$. 
Let $f$ be diagonal-free. Then if $j_1,\dots$, $j_k \le r$, then
\begin{equation}\label{inter}
\E[Q(f,\,\bm\xi_{j_1},\ldots,\bm\xi_{j_k})\,|\,\ca G_r]=
r^{-k}Q(f,(\bm\xi_1+\ldots+\bm\xi_r)^{k}).
\end{equation}
We should point out that the last term represents the random polynomial
$$ Q(f,(\bm\xi_1+\ldots+\bm\xi_r)^{k}) =
   \sum_\ik f_\ik (\xi_{1 i_1} + \dots + \xi_{r i_1}) \cdots
   (\xi_{1 i_k} + \dots + \xi_{r i_k}) .$$
Equation \refrm{inter} follows because 
\[
\E  (f_\ik \xi_{j_1 i_1} \cdots \xi_{j_k i_k} | G_r)
= f_\ik \E(\xi_{j_1 i_1} | G_r) \cdots\E(\xi_{j_k i_k} | G_r)
\]
because $f$\ is diagonal free, and hence $\xi_{j_1 i_1},\dots$, $\xi_{j_k i_k}$
are independent if $f_\ik \ne 0$, and also because for $j \le r$
\[
\E(\xi_{j i} | G_r) = r^{-1}(\xi_{1 i} + \dots + \xi_{r i}) .
\]
We also point out the following easy consequence of the triangle inequality
for $\L$.
\begin{equation}\label{(0)}
\lnorm \wh{Q}(f,\BX) \rnorm_{\L(\BE)} \le \lnorm Q( f,\BX)\rnorm_{\L(\BE)}.
\end{equation}

Now we are ready to present our first decoupling inequality.  This result
allows us to decouple random polynomials in the rearrangement invariant
norm.
\begin{theorem}\label{AB} 
Let $f = (f_\ik)$ be a diagonal free array of vectors from $\BE$. 
Let $\bm\xi, \bm\xi_1,\ldots,\bm\xi_k$
be sequences of integrable independent real random variables.  
Let $\L$ be a r.i. space of random variables, containing 
$\xi_{1}\cdots \xi_{k}$ (hence, norms of all finitely
supported polynomials spanned by
 $\bm\xi_1,\ldots,\bm\xi_k$).
\begin{itemize}
\item[{\rm (A)}] 
Assume that $\bm\xi,\bm\xi_1,\bm\xi_2,\ldots$ are independent and
identically distributed. 
Then
\[
\lnorm   Q(f,\bm\xi^{k})\rnorm _{\L(\BE)}
\le A\lnorm  Q(f,\bm\xi_1,\ldots,\bm\xi_k)\rnorm  _{\L(\BE)}
\]
where $A=A_k\sim (2k)^k$, or, if $\E\bm\xi=0,\, A_k=k^k$.
\item[{\rm (B)}] 
Assume that $\bm\xi,\bm\xi_1,\ldots,\bm\xi_k$ are interchangeable. 
Then
\[
\lnorm  \wh{Q}(f,\bm\xi_1,\ldots,\bm\xi_k)\rnorm_{\L(\BE)} \le
B\lnorm {Q}(f,\bm\xi^{k})\rnorm_{\L(\BE)} ,
\]
where $B=B_k\sim k^k/k!$.
\end{itemize}
\end{theorem}
\Proof

During this proof, we will suppress the subscript $\L(\BE)$ on the norms.

{\bf (A)}: 

\underline{ Step 1$^o$: Centering procedure}: 
\vspace{5pt}

\def\ir{{i_1,\dots,i_r}}

Denote $\ov{\xi}=\xi-\E[\xi]$,  $\bm m=(m_1,m_2,\ldots)$,
where $m_i=\E[\xi_i]$. 
For $1\le r\le k$, if $f = (f_\ir)$, then we have
\begin{equation}\label{centdec}
\lnorm  Q(f,\ov{\bm\xi}_1,\ldots,\ov{\bm\xi}_r)\rnorm 
\le 2^r\,\lnorm Q(f,\bm\xi_1,\ldots,\bm\xi_r)\rnorm .
\end{equation}

Indeed, by interchangeability
\[
\begin{array}{c}
\displaystyle
\lnorm  Q(f,\ov{\bm\xi}_1,\ldots,\ov{\bm\xi}_r)\rnorm =
\lnorm  
Q(f,\bm\xi_1-\bm m,\ldots,\bm\xi_r-\bm m)\rnorm \\
\displaystyle
=\lnorm  \Sum_{(\delta_1,\dots,\delta_r)\in\set{0,1}^r}
\sum_\ir f_\ir \xi_{1 i_1}^{\delta_1} \cdots \xi_{r i_r}^{\delta_r}
m_{1 i_1}^{1-\delta_1} \cdots m_{r i_r}^{1-\delta_r}\rnorm \\
\displaystyle
\le \sum_{j=0}^r {r \choose j} \lnorm \sum_\ir f_\ir
\xi_{1 i_1} \cdots \xi_{j i_j} m_{j+1, i_{j+1}} \cdots m_{r i_r} \rnorm \\
\displaystyle
= {\Sum_{j=0}^r}{r \choose j}
\lnorm  Q(f, \bm\xi_1,\ldots,\bm\xi_j, \bm m^{(r-j)})\rnorm .\\
\end{array}
\]

The latter expression is equal to 
\[
{\Sum_{j=0}^r}{r \choose j}
\lnorm  Q( f, \bm\xi_1,\ldots,\bm\xi_j, 
\E[\bm\xi_{j+1}|\BX_j],\ldots,
\E[\bm\xi_r|\BX_j])\rnorm
={\Sum_{j=0}^r}{r \choose j}
\lnorm  \E[Q( f, \bm\xi_1,\ldots,\bm\xi_j, 
\bm\xi_{j+1},\ldots,
\bm\xi_r) | \BX_j ]\rnorm .
\]
where $\BX_j$ is the matrix $[\bm\xi_1,\ldots,\bm\xi_j]$. 
Using the contractivity property \refrm{L}, we
estimate the above term from above by
\[
{\Sum_{j=0}^r}{r \choose j}
\lnorm  Q(f,\bm\xi_1,\ldots,\bm\xi_r)\rnorm 
= 2^r \lnorm  Q(f,\bm\xi_1,\ldots,\bm\xi_r)\rnorm.
\]
\vspace{5pt}

\underline{Step 2$^o$: proving (A)}: 
Arguing similarly to Step 1$^o$ above, we note that
\begin{eqnarray*}
\lnorm  Q(f, \bm\xi^k)\rnorm =
\lnorm   Q(f, (\ov{\bm\xi}+\bm m)^k)\rnorm 
\le \Sum_{r=0}^k {k \choose r} \lnorm  
Q( f, \ov{\bm\xi}^r, \bm m^{k-r})\rnorm := Q_0.
\end{eqnarray*}

Now, using \refrm{L}, and noting that 
$\E(\ov{\bm\xi}_1 + \cdots + \ov{\bm\xi}_r | \ov{\bm\xi}_1) = \ov{\bm\xi}_1$, 
which is identically
distributed to $\ov{\bm\xi}$, it follows that

\begin{eqnarray*}
\displaystyle
Q_0 \le  \Sum_{r=0}^k {k \choose r} \lnorm   
Q(f, (\ov{\bm\xi}_1+\cdots +\ov{\bm\xi}_r)^r,
\bm m^{k-r})\rnorm 
\end{eqnarray*}
Then, in virtue of \refrm{inter}, the latter expression is equal to 
\[
\Sum_{r=0}^k {k \choose r} \lnorm   \,r^r\,
\E[\,Q(f,\ov{\bm\xi_1},\ldots, \ov{\bm\xi_r}, 
\bm m^{k-r})\,|\,\ca G_r]\rnorm 
=:Q_1
\]

Using \refrm{L}  and  applying the centering procedure
\refrm{centdec}, the above term gets the following upper bounds:

\begin{eqnarray*}
Q_1\le 
\Sum_{r=0}^k {k \choose r} \lnorm  
 \,r^r\,
Q( f,\ov{\bm\xi_1},\ldots, \ov{\bm\xi_r}, 
\bm m^{k-r})\rnorm \\
\displaystyle
\le \Sum_{r=0}^k {k \choose r}  (2r)^r
\lnorm  \, \,
Q(f,\bm\xi_1,\ldots, \bm\xi_r, \bm m^{k-r})\rnorm 
=: Q_2
\end{eqnarray*}

Then, by interchangeability, independence of columns,  and \refrm{L} again, we
keep estimating, as follows
\begin{eqnarray*}
Q_2= \Sum_{r=0}^k {k \choose r}   (2r)^r\,
\lnorm  
\E[\,Q(f,\bm\xi_1,\ldots, \bm\xi_r, 
\bm\xi_{r+1},\ldots,\bm\xi_k)\,|\,\BX_r\,]\rnorm \\
\le \Sum_{i=0}^k {k \choose r} (2r)^r\, 
\lnorm  
Q( f,\bm\xi_1,\ldots, \bm\xi_r, 
 \,\bm\xi_{r+1},\ldots,\bm\xi_k)\rnorm \\
=A_k \lnorm  \,
Q( f,\bm\xi_1,\ldots,\bm\xi_k)\rnorm .\\
\end{eqnarray*}
Clearly, $A_k\le (2k+1)^k$ (notice that $A_k\ge c (2k)^k$). If $\E\xi=0$, the 
use of centering procedure, and the triangle inequality is superfluous, 
hence the constant decreases to $k^k$.
\vspace{5pt}

(B): By the Mazur-Orlicz polarization formula \refrm{MO}, and \refrm{inter}, 
we obtain the following bounds
\[
\begin{array}{rl}
\lnorm  \wh{Q}(f,\bm\xi_1,\ldots,\bm\xi_k)\rnorm 
=&\lnorm  \Rec{k!} Q(f, \Sum_{\ddl} 
(-1)^{k-|\ddl|}
(\dl_1\bm\xi_1+\ldots+\dl_k\bm\xi_k)^k)\rnorm \\
\le &\Sum_{\ddl}\lnorm  {\Rec{k!}} 
Q( f, (\dl_1\bm\xi_1+\ldots+\dl_k\bm\xi_k)^k )\rnorm \\
=&\displaystyle 
\Sum_{r=0}^k{k \choose r}
\lnorm  \Rec{k!}Q(f, (\bm\xi_1+\ldots+\bm\xi_r)^k)\rnorm \\
\le &\displaystyle
\Sum_{r=0}^k{k \choose r}
\lnorm  \Frac{r^k}{k!} Q(f,\bm\xi_1^k)\rnorm \\
= &
 B_k\lnorm Q(f,\bm\xi_1^k)\rnorm. \\
\end{array}
\]
That $B_k\sim k^k/k!$ is easy to verify.
The proof is completed. 
\QED

\subsection{Extended multilinear forms and U-Statistics}

In this section, we show how to extend Theorem \ref{AB} to the so called
U-Statistics.  Let $k$\ be a positive integer.  Let $F = (F_\ik)$\
be an array of strongly Borel measurable functions $F_\ik:\R^k \to \BE$
such that
\[
\begin{array}{rl}
\mbox{\rm \hspace{-10pt}  (F0)\hspace{10pt} }&
           \mbox{\sf $F_\ik = 0$ if some $i_j$\ 
and $i_{j'}$ are identical for $j \ne j'$}\\
\mbox{\rm \hspace{-10pt}  (F1)\hspace{10pt} }&   
           \mbox{\sf $F_\ik=0 $ for all but finitely many $(\ik)$;}\\
\end{array}
\]
Then we are going to consider U-Statistics, that is,
expressions of the form
\[
F(\bm\xi_1,\dots,\bm\xi_k)\df\sum_{\ik} F_\ik(\xi_{1i_1},\ldots,\xi_{ki_k}),
\]
where $\bm\xi_1,\dots,\bm\xi_k$ are real valued random variables.  (Here, $\R$\
could be replaced with any other measure space, but there is no loss of
generality to take it as $\R$.)  We are going to exercise the same notational
devices as for the random polynomials, so that the undecoupled U-Statistic is
written
\[
F(\bm\xi^k)\df\sum_{\ik} F_\ik(\xi_{i_1},\ldots,\xi_{i_k}) .
\]
As before, in order to prove the results, we require certain symmetry
properties to hold for $F$.  So we defined the symmetrized version of
$F$ as follows:
\[
\hat F(x_1,\dots,x_k)
\df {1\over k!} \sum_\sigma
F_{i_{\sigma_1},\dots,i_{\sigma_k}}(x_{i_{\sigma_1}},\ldots,xi_{i_{\sigma_k}}),
\]
where the sum runs over all permutations of $[1,k]$, and we set
\[
\hat F(\bm\xi_1,\dots,\bm\xi_k)
\df \sum_{\ik} 
\hat F_\ik(\xi_{1i_1},\ldots,\xi_{ki_k}) .
\]
Decoupling results were proved by \cite{Pen} for Orlicz modulars (and
so by Note~\ref{riorl} below, one can obtain 
results for all rearrangement invariant
spaces).  We will prove similar decoupling results, weakening some of
the hypotheses.  

More interestingly, we are going to prove the decoupling results for 
U-statistics as a corollary of Theorem~\ref{AB} which decouples random
polynomials.  The technique is to approximate the U-statistic as
a sum of random polynomials.  That is, let $D$\ be an integer,
and for $1 \le d \le D$, let $f^d = (f^d_\ik)$ be a diagonal free array
of vectors in $\BE$\ taking only finitely many non-zero values, and let
$(\bm\xi^d_1:1\le d \le D),\dots$, $(\bm\xi^d_k:1\le d \le D)$ be 
sequences of independent random
variables.  Then we set
\begin{equation}\label{newtensor}
R(f;\bm\xi_1,\dots,\bm\xi_k) \df
\sum_{d=1}^D Q(f^d;\bm\xi^d_1,\dots,\bm\xi^d_k) .
\end{equation}

Then the remarkable thing is that the proof of Theorem~\ref{AB} works for
$R(f;\bm\xi_1,\dots,\bm\xi_k)$\ exactly as it does for 
$Q(f;\bm\xi_1,\dots,\bm\xi_k)$, that is, we have the following result.

\begin{theorem}\label{ABnew}
Theorem \ref{AB} is valid, for the
multilinear form \refrm{newtensor}.
\end{theorem}

A version of the following result for Orlicz modulars
$\E\phi(\cdot)$, where $\phi$ was a moderately increasing function,
was proved in \cite{McCT:bandec}. In that paper, terms of the underlying sums
were sign-randomized, i.e., each $F(\bi, \cdot)$ was multiplied by  Walsh
functions $\eps_{i_1}\cdots \eps_{i_k}$. More precisely, the decoupling was
proved for
\begin{equation}\label{Fbeps}
(F\circ\beps)(\bm\xi_1,\dots,\bm\xi_k)\df
\sum_{\bi} \eps_{i_1}\cdots\eps_{i_k}F_\ik(\xi_{1i_1},\ldots,\xi_{ki_k}),
\end{equation}
where $\beps$ is independent of $\BX$.  That the
presence of Walsh functions is not necessary in the context of Orlicz modulars,
was shown in \cite{Pen}. We observe that the following result, generalizing
theorems in the mentioned papers, is implicit in the main decoupling principle.
Moreover,  constants remain the same. For the sake of completeness we give the
full proof.

In the proof we will use the fact that 
any inequality involving norms of functions of discrete 
r.v.'s, that converge to some limits, is preserved for these 
limits. 
Fix $\bm i$, say $\bm i=(1,\ldots,k)$. Consider $X=F(\xi_1,\cdots,\xi_k)$. 
We may assume that the probability space $(\Omega,\ca F,\P)$  and $\ca F$ is 
spanned by $\xi_1,\xi_2,\ldots$. Also, we may assume that it is separable, 
i.e., $\ca F_n=\sigma\{\bigcup_n \ca F_n\}$, where $\ca F_n$ are finite 
$\sigma$-fields.  Put $\xi_i^n=\E[\xi|\ca F_n]$. Thus 
$\E[X|\ca F_n]\to \E[X|\ca F]=X$ a.s. and in $\L$.

\begin{theorem}\label{F}
Let $F:\N^k\times\R^k\to\BE$  satisfy {\rm (F0) -- (F1)}, and also
and additional condition:
\[
\begin{array}{rl}
\mbox{\rm \hspace{-10pt}  (F2)\hspace{10pt} }&
           \mbox{\em $F(\bi;\xi_{i_1},\ldots,\xi_{i_k})\in \L$ 
                                      for every $\bi\in\N^k$.}\\
\end{array}
\]
Let
$\bm\xi_1,\bm\xi_2,\ldots$ be sequences of independent
random variables.  
\begin{itemize}
\item[(A$'$)] 
Let $\|\cdot\|$ be a r.i. norm. Let $\bm\xi_1,\bm\xi_2,\ldots $ be 
independent and identically
distributed. Then
\[
\left\|\,F(\bm\xi^k)\,\right\|\le A 
\left\|\,\,F(\BX)\,\right\|,
\]
where $A$ is the constant from Theorem \ref{AB}.(A).
\item[(B$'$)] 
Let $\bm\xi_1,\ldots,\bm\xi_k$ be interchangeable (in particular, i.i.d.)
and $\|\cdot\|$ be a r.i.\ norm.  Then
\[
\left\lnorm \,\wh{F}(\BX)\,\right \rnorm
\le B_k \left\lnorm \,F(\bm\xi^k)\,\right\rnorm.
\]
where $B$ is the constant from Theorem \ref{AB}.(B).
\end{itemize}
\end{theorem}
\Proof 
By Note~\ref{riorl}, 
we may assume that the rearrangement space $\L$\ is separable.
In that case, we may assume without loss of generality that $\bm\xi_1,
\bm\xi_2,\ldots$\ are real discrete random variable.

Thus
we may assume that the random variables are defined on a product probability
space
\[
(\prod _{ij}\Omega_{ij}, \ca
(F^{\otimes N})^{\otk},(\P^{\otimes \N})^{\otk}), 
\]
where $\Omega_{ij}$ are
equal, $\Omega_i=\prod_j\Omega_{ij}$, and the superscript $\otimes$
indicates the product $\sigma$-field and the product probability,
respectively. So, let 
\[
\xi_{i}=\sum_m x_{im}\I{A_{im}},
\]
where $A_{i1},A_{i2},\ldots\subset\Omega_i$ are bases of rectangular sets that 
form a disjoint finite partition of $\Omega$, and let $(A_{sim})$,
$s=1,\ldots,k$, be independent copies of $(A_{im})$.  Put
$I_{sim}=\I{A_{sim}}$, hence
\[
\xi_{si}=\sum_m x_{im} I_{sim},\qquad s=1,\ldots,k.
\]
Then
\[
\sum_{\ik}F_{ik}(\xi_{1i_1},\ldots,\xi_{ki_k})=
\sum_{m_1,\dots,m_k} \sum_{\ik}
F_\ik(x_{i_1m_1},\ldots, x_{i_km_k})
I_{1i_1m_1}\cdots I_{ki_km_k}.
\]
Now, we can apply Theorem \ref{ABnew}, and the proof is complete.
\QED

\subsection{An example in a certain Lorentz space}\label{lorentzeg}

 Motivated by the results in \cite{Pen:Wald} where
the problem as to when expectation results imply tail probability results is
treated, 
we obtained the following asymptotic tail probability comparison.

\begin{proposition}\label{limsup}%
Let $\xi$ or $\eta$ be the norm 
of $\wh{F}(\bm\xi^k)$ or $\wh{F}(\BX)$, and let $W:[0,\infty) \to [0,\infty)$\
be an increasing function
such that there exists constants $p>1$\ and $c>0$\ such that
$W(st) \le c s^p W(t)$\ for all $0<s<1$\ and all $t>0$.  Then there is
a constant $C$, depending only on $p$\ and $c$, such that
\[
\limsup_{t\to \infty} W(t) \P(\xi \ge t) \le 
\limsup_{t\to\infty} W(t) \P(\eta \ge C t ). 
\]
\end{proposition}

This result is a consequence of Theorem~\ref{F}, 
and follows by arguments from the theory of Lorentz-Zygmund spaces. 
If $\xi$ is a random variable, let $F(t)=\P(|\xi|\ge t)$, and define
the decreasing
rearrangement of $\xi$\ to be the
function $\xi^*(t)\df \sup\{s:F(s) > t \}$\ (i.e. the
right-continuous inverse of $F$). 
Obviously, $|\xi|$ and $\xi^*$ are equidistributed. 
When $\xi$ is integrable,  an average operator is often considered
\[
\xi^{**}(t)=\frac{1}{t}\int_0^t \xi^*(u)\, du,
\]
which corresponds to  a rearrangement invariant norm for every $t>0$. 
Therefore, our decoupling inequalities for U-statistics hold for 
$\Phi(X)=(\|X\|)^{**}$. That is, denoting by $\xi$ or $\eta$ the norm 
of $\wh{F}(\bm\xi^k)$ or $\wh{F}(\BX)$, we have
\[
\xi^{**}(t) \le C \eta^{**}(t)
\]
for some constant $C>0$. 

Now consider the Lorentz-Zygmund space  defined by the quasi-norm
$|||f|||=\sup_x w(x) \xi^*(x)$, where $w:[0,1]\to 
[0,\infty)$\ is an increasing function.  
Note that $|||f||| \le 1$\ if and only if
$ \sup_t W(t) \P(|f| \ge t) \le 1$, where $w(t)=\Frac{1}{W^{-1}(1/t)}$.
If $W$\ satisfies the relation given in Proposition~\ref{limsup}, then 
for some constant $c$, the function $w$\ satisfies
the relation $w(x)\le c a^{-1/p} w(xa)$\ for $a \le 1$.  Then it is possible
to show that $|||f||| \le |||f^{**}||| \le C |||f|||$.  Indeed, the first
inequality is obvious, and for the second:
\[
w(x) f^{**}(x) = w(x) \int_0^1 f^*(xa) \, da
\le c \int_0^1 a^{-1/p} w(xa) f^*(xa) \, da
\le \frac {cp}{p-1} ||| f ||| .
\]
Thus, to show Proposition~\ref{limsup}, 
let $F(t)=\P(\xi\ge t)$\ and $G(t)=\P(\eta\ge t)$.
Then
\[
\sup_t W(t) F(t)\le 
\sup_t W(t) G(Ct) .
\]
If we now set $w(x)=0$, for $x\ge x_0$, the same argument applies, 
and letting $x_0 \to \infty$, we obtain
\[
\limsup_{t\to \infty} W(t) F(t)\le \limsup_{t\to\infty} W(t) G( C t ). 
\]


\section{A discourse on probability tails}
\subsection{$L^p$-estimates imply tail estimates}
\subsubsection{Auxiliary results}
The following result can be found in \cite{AsmM}.
\begin{lemma}\label{Lemma 2.1}   
Let $\{X;~X_i\}$ be a sequence of positive i.i.d.\ random variables.
Then, for all positive integers $n$,
all $\alpha > 0 $
and all $0\le \theta \le n$
\begin{eqnarray*}
\P(X \ge \alpha) \ge {\theta \over n} \quad 
&\Rightarrow \quad \P(\sup_{1\le i\le n} X_i  \ge \alpha) 
\ge {\theta\over 1+\theta},&  \\
\P(X \ge \alpha) \le {\theta \over n} \quad
&\Rightarrow \quad
\P(\sup_{1\le i\le n} X_i \ge \alpha )  \le \theta. &   
\end{eqnarray*}
\end{lemma}                                         

\Proof
To show the first inequality, observe first that for $\theta > 0$,
$$ (1 - {\theta \over n} )^n \le {1 \over (1+\theta ) }.$$
Hence, by independence assumption, 
\[
\begin{array}{rl}
 \P(\sup_j X_j \ge \alpha ) 
      &= 1 - \P(\sup_j X_j < \alpha ) \\
      &= 1- \displaystyle\prod_{j=1}^n \P(X_j < \alpha ) r\\ 
      &\displaystyle\ge 1- (1-{\theta \over n} )^n  
\ge 1-{1\over (1+\theta )} = {\theta \over (1+ \theta )}. \\
\end{array}
\]
The second inequality is easy: from the imposed condition, one gets
$$
\P(\sup_{1\le j\le n}X_j \ge \alpha ) \le 
\sum_{i=1}^n \P( X_i  \ge \alpha ) \le \theta .
$$
The proof is completed.
\QED

The following result can be found in \cite[Chap. 4]{LedT}.

\begin{lemma}\label{Lemma 2.2}
Consider a positive random variable $Z $ such that
$ \|Z\|_q \le C\|Z\|_p$\ for $q>p>0$. Then,
\[
\P(Z>t) \le (2C^p)^ {q\over p-q} \qquad
\Rightarrow \qquad \|Z\|_p \le 2^{1/p} t 
\quad\hbox{and}\quad
\|Z\|_q \le 2^{1/p} Ct.
\]
\end{lemma}
Putting together Lemmas \ref{Lemma 2.1} and \ref{Lemma 2.2} we get the
following.

\begin{lemma}\label{Lemma 2.3}%
Let $\{X;~X_i\}$ be a sequence of positive i.i.d.\ random variables.
Assume that there exists a constant $c$ such that for $0<p<q<\infty$ 
$$
\|\sup_{1\le i \le n}\| X_i\|\|_q \le c\|\sup_{1\le i \le n}\|X_i\|\|_p.
$$ 
Then, letting
$\theta = (2c^p)^{q\over p-q} $, 
\[
\P(X\ge t) \le 
{\theta \over n} \quad\Rightarrow\quad
\|\sup_{ 1\le i\le n} X \|_p \le 2^{1\over p}t. 
\]
\end{lemma}                                         

For later reference, we also include the next lemma.
\begin{lemma}\label{Lemma 2.4}                      
Let $\{X;~X_i\}$ be a sequence of positive i.i.d.\ random variables.
Then,
\[
\|\sup_{1\le i \le n} X_i \|_p
\le t
\quad\Rightarrow\quad \P(X\ge 2^{1\over p} t) \le {1 \over n}.
\]
\end{lemma}                                         

\Proof 
Use Chebychev's inequality and Lemma \ref{Lemma 2.1} with $\theta = 1$.
\QED

\subsubsection{Main result}
Now we are ready to prove an  extension of a result from \cite{AsmM} 
that  deals
with strict tail probability comparisons for
pairs of random variables.

\begin{theorem}\label{2.1}
Let $(X,\,X_i)$ and $(Y,\, Y_i)$ be sequences of positive
i.i.d.\ random variables. 
For some $0<p<q$\ and all positive integers $n$ assume that

\begin{equation}\label{(0.1)}
\|\sup_{1\le i \le n} X_i\|_q \le c_1\|\sup_{1\le i \le n} X_i \|_p, 
\end{equation}
and 
\[
\|\sup_{1\le i \le n} Y_i\|_p \le c_2\|\sup_{1\le i \le n} X_i\|_p.
\]
Then there exists $c_3$, depending only on $p$, $q$, $c_1$\ and $c_2$\  
such that for all $t \ge 0$
\[
\P(Y\ge c_3t) \le c_3 \P(X\ge t).
\]
\end{theorem}

\Proof 
Given an arbitrary $\alpha = \alpha_1 > 0$ with 
\begin{equation}\label{(2.6)}
\P(Y\ge \alpha_1) > 0, 
\end{equation}
choose $\mu $ to be the 
smallest positive integer satisfying
\begin{equation}\label{(2.7)}
{1\over 2\mu} \le \P(Y\ge \alpha_1)\le {1\over \mu}. 
\end{equation}
>From Lemma \ref{Lemma 2.1} it follows that
$$
\P(\sup_{1\le j \le \mu}Y_j\ge \alpha_1) \ge {1\over 3}.
$$ 
Hence, by Chebychev's inequality,
$$
{1\over \alpha_1^p }
\|\sup_{1\le j \le \mu}Y_j\|_p^p \ge {1\over 3}, 
$$
which, by assumptions, yields
$$
{c_{2}^p\over \alpha_1^p }
\|\sup_{1\le j \le \mu}X_j\|_p^p \ge {1\over 3}, 
$$
and, consequently, for any $\alpha_2 > 0$, 

$$ 
{1\over \alpha_2^p }\|\sup_{1\le j \le \mu}X_j\|_p^p \ge 
{\alpha_1^p\over 3\alpha_2^p c_{2}^p}. 
$$
In particular, if $\alpha_2^p = {\alpha_1^p\over 6 c_{2}^p}$, 
we get from the latter inequality that 
$$
{1\over \alpha_2^p }
\|\sup_{1\le j \le \mu}X_j\|_p^p \ge 2. 
$$
Now, Lemma 2.3 implies that
$$
\P(X\ge {\alpha_1\over (6^{1\over p} 
c_{2} )}) = \P(X\ge \alpha_2)\ge  
 (2 c_1^p)^{q\over p-q} {1\over \mu} . 
$$ 
Finally, \refrm{(2.7)} gives,
\begin{equation}\label{(2.14)}
\P(X\ge {\alpha_1\over (6^{1\over p} c_{2})}) \ge
 (2 c_1^p)^{q\over p-q} \P(Y \ge \alpha_1 ). 
\end{equation}
Note that \refrm{(2.14)} holds for all $\alpha_1$ for which \refrm{(2.6)}
holds. For any other $\alpha_1 > 0$, \refrm{(2.14)} holds trivially.
\QED

\subsubsection{Contraction for multipliers}

Condition \refrm{(0.1)} yields an example of a class of random variables
with the so called Marcinkiewicz-Paley-Zygmund property (MPZ in short). The
concept  was studied in \cite{KraS:hyp}, and can be traced back to  \cite{PalZ}
and \cite{MarZ}. A family $\ca Z\in L^q_+$ of random variables is said to be in
the class $MPZ(q)$ (in short: {\em have MPZ}), if one of the following 
equivalent conditions is satisfied:

\begin{eqnarray}
&\exists~ p<q \mbox{ (equivalently, $\forall~ q\le p$)~ }~
m_{q,p}\df\sup_{Z\in \ca Z} \Frac{\|Z\|_q}{\|Z\|_p}<\infty
\label{MPZ1}&\\
\rule{6mm}{0pt} \nonumber\\
&\exists ~ \delta>0~~
\inf_{Z\in \ca Z} \P(Z>\delta\|Z\|_q)>\delta&
\label{MPZ2}
\end{eqnarray}

That is, \refrm{(0.1)} involves $\ca Z=\set{\sup_{1\le i\le n} \|X_i\|:n\in
\N}$. Also, in \cite{KraS:hyp} it was shown that the space of diagonal-free
random polynomials of finite degree, spanned by symmetric random variables
with, so called, semi-regular distribution, 
has MPZ. A random variable $\xi$ is
said to have the {\em semi-regular distribution}, if its tail 
$G(t)=\P(|\xi|>t)$
satisfies the relation,
\[
V(a)=\limsup_{t\to\infty} G(at)/G(t)<1
\]
for some (or all) $a>1$ (by convention, $0/0=0$). For example, any bounded
random variable has semi-regular distribution. In particular, (the norm of) any
normed space-valued Rademacher polynomial of degree $d$ has MPZ with the
constant $m_{qp}=[2 (q-1)/(p-1)]^{d}$ (\cite[Corollary 2.7]{KraS:hyp}).

Now, the essence of Theorem \ref{2.1} is that  the continuity of a certain
operator, once is fulfilled by means of $L^p$-norms, will be also fulfilled by
means of probability tails. We will illustrate this concept by the following
result.
\begin{theorem}\label{contrRad}
Let $f_\ik$ be a, 
finitely supported, diagonal-free
array, taking values in a Banach space
$\BE$.  Let $\bm\xi = (\xi_1,\dots)$\ be a sequence 
of symmetric independent random variables. 
\begin{enumerate}
\item[{\rm (i)}] ({\bf Contraction Inequality})
There is a constant $c>0$ such that 
\[
\P(\|Q(f, \bm (s\bm \xi)^k)\|>c t)\le 
c \P(\|Q(f, \bm \xi)\|>t),\qquad t>0,
\]
where $\bm s\bm\xi=(s_i\xi_i)$ and $\|\bm s\|_\infty=\sup_i |s_i|\le 1$. 
\item[{\rm (ii)}] ({\bf Maximal Inequality})
There is a constant $C>0$\ such that
\[
\P(\sup_{m_1,\dots,m_k }\|T_{m_1,\dots,m_k}Q(f, \bm \xi^k )\|>C t) \le C
\P(\| Q( f, \bm \xi^k)\|>t),\qquad t>0,
\]
where 
\[
T_{m_1,\dots,m_k} Q(f;\BX) = \sum_{i_1\le m_1,\dots, i_k\le m_k} 
f_{\ik} \xi_{1i_1}\cdots\xi_{ki_k}.
\]
\end{enumerate}
\end{theorem}
\Proof
Let us first prove~(i).
In virtue of symmetry assumption and Fubini's theorem, it suffices to 
give the proof for the 
case when $\BX$\ is a matrix of Rademacher random variables.  
Let $Q_1,\dots,Q_n$\ be independent copies of $Q(f,\bm)$.
Then the vector $(Q_1,\dots,Q_n)$\ is a Rademacher homogeneous
polynomial of degree $k$\ taking values in $\ell_n^\infty(\BE)$.
Similarly, let $R_1,\dots,R_n$\ be independent copies of $Q(f, \bm s\bm \xi)$.
>From the contraction principle for $L^p$-norms, which may be found in
\cite[Remark 2.9]{KraS:hyp} (essentially, it is due
to \cite{Kwa:dec}), it follows that for all $p \ge 1$
\[
\| \|(R_1,\dots,R_n)\|_{\ell_n^\infty(\BE)} \|_p
\le c
\| \|(Q_1,\dots,Q_n)\|_{\ell_n^\infty(\BE)} \|_p .
\]
>From the observation that $\|(Q_1,\dots,Q_n)\|_{\ell_n^\infty(\BE)}
= \sup_{1 \le i \le n} \|Q_i\|_{\BE}$, and using the fact that 
Rademacher polynomials are MPZ, and also citing Theorem~\ref{contrRad} 
above, the result
follows.

The proof of part~(ii) is the same, using the corresponding result
for $L^p$-norms of polynomials for the Rademacher random variables, 
which follows easily from
\cite{McCT:realdec} and L\'evy's inequality.

\begin{remark}\label{contrRAD:rem}\rm
In fact, Theorem~\ref{contrRad} is also valid for the sign randomized 
U-Statistics as in equation~\refrm{Fbeps}.  The proof is identical.
\end{remark}

\begin{theorem}\label{compRad}[ Comparison Inequality]
Let $(f_{i_1,\dots,i_k})$ be a diagonal-free, finitely supported, diagonal free
array.  Let 
$\bm\xi = (\xi_i)$\ and $\bm\eta = (\eta_i)$\ be sequences
of symmetric independent random variables. 
such
that, for some constant $A>0$,
\[
\P(|\xi_i|> t)\le A\P(|\eta_i|>t),\qquad t>0, ~i\in \N.
\]
Then, for some constant $K=K(c,d,A)$,
\[
\P(\|Q(f, \bm \xi^k)\|> t)\le
K \P(K \|Q(f, \bm \eta^k)\|>t),\qquad t>0.
\]
\end{theorem}
\Proof
We have 
\[
Q(f, \bm \xi^k)\D Q(f, (\beps|\bm \xi|)^k),
\]
and 
\[
Q(f, \bm \eta^k)\D Q(f, (\beps|\bm \eta|)^k),
\]
where $\beps$ is a Rademacher sequence independent of $\bm\xi$ and $\bm\eta$. 

If $A=1$, then we may replace each $|\xi_i|$ and $|\eta_i|$ by their decreasing
rearrangements $|\xi_i|^*$ and $|\eta_i|^*$, respectively. The assumption yields
$|\xi_i|^*\le |\eta_i|^*$ a.s. Hence, by Theorem~\ref{contrRad}(i), 
the inequality follows.

Let $A>1$. Then there exist a sequence $\bm\alpha=(\alpha_i)$ of i.i.d. random
variables, independent of $\bm\xi$, such that
$\P(\alpha_i=1)=1/K,\,\P(\alpha_i=0)=1-1/K$, so that
$\P(\alpha_i|\xi_i|>t)=\P(|\xi_i|>t)/K$. Therefore, by the first part of the
proof,
\[
c\P(\|Q(f, (\beps\bm\alpha|\bm \xi|)^k)\|> t)\le
\P(K \|Q(f, (\beps|\bm \eta|)^k)\|>t),\qquad t>0.
\]
Conditioning on $\bm\xi$, it remains to prove that, 
for every polynomial $Q$, and every diagonal free array $f_\ik$,
\begin{equation}\label{compeps}
\P(\|Q(f, \beps^k)\|> m t)\le 
m\P(\|Q(f, (\beps\bm\alpha)^k)\|> t)
\end{equation}
for some constant $m=m_k$.
Let $\bm\beta=\beps\bm\alpha$.  Then it is clear that $\bm\beta$\ is 
semiregular, as defined earlier, and hence homogeneous random 
polynomials of degree $k$\
over $\beta$\ have MPZ.  Furthermore, the comparison inequality is
true for $L^p$\ for $p\ge 1$\ (see, for example,
\cite[Theorem~2.13]{KraS:hyp})  Hence arguing as in the proof of
Theorem~\ref{contrRad}, we obtain \refrm{compeps}.
\QED

\subsection{Decoupling for tails}
In order
to prove any tail inequality of the type $\P(\xi>t)\le K\P(\eta>t)$, where
$\xi,\eta$ are real random variables, it is enough to prove it for an
arbitrarily chosen conditional probability
\[
\P \left[\,\xi>t\,|\,\ca G \,\right]\le 
K\P \left[\,\eta>t\,|\,\ca G\,\right].
\]
This observation was used in proving the inequality 
\cite[(6.9.5)]{KwaW:book}.
Denote by $\ca G$ the
$\sigma$-field spanned by all random variables of the form $\sum_{j=1}^i
h(\bm\xi_j)$ (in other words, by the random point measure
\[
 \sum_{j=1}^i\delta_{\bm\xi_j}
\]
on $(\R^\N)^k$ (cf. \cite[p.\ 182]{KwaW:book}). Then 
$(\bm\xi_1,\ldots,\bm\xi_k)$ is concentrated on a finite 
permutation invariant subset of  $(\R^\N)^k$. Now, \refrm{inter} 
can be rewritten, as follows (recall the notation, preceding \refrm{inter}).
\begin{equation}\label{newinter}
\E[Q(f,\bm\xi_{j_1},\ldots,\bm\xi_{j_k})\,|\,\ca G]=
k^{-k}Q(f,(\bm\xi_1+\ldots+\bm\xi_k)^{k}),
\end{equation}

\begin{theorem}
\label{ABtail} 
Let $f,\bm\xi,\bm\xi_1,\ldots,\bm\xi_k$ be as in 
Theorem \ref{AB} (but \underline{we do not assume integrability}). 
\begin{itemize}
\item[{\rm (A$''$)}] 
Let $\bm\xi,\bm\xi_1,\bm\xi_2,\ldots$ be 
independent and symmetric. Then, there exists a constant $A''$,
depending only on $k$, such that, for all $t\ge 0$,
\[
\P(\lnorm   Q( f,\bm\xi^{k})\rnorm\|\ge A'' t)\le 
A''\P(\lnorm  Q( f,\bm\xi_1,\ldots,\bm\xi_k)\rnorm\ge t). 
\]
\item[{\rm (B$''$)}] 
Let $\bm\xi,\bm\xi_1,\ldots,\bm\xi_k$ be interchangeable. 
Then, there exists some constant $B''$, 
depending only on $k$, such that, for all $t\ge 0$,
\[
\P(\lnorm  \wh{Q}(f,\bm\xi_1,\ldots,\bm\xi_k)\rnorm\ge B''t) 
\le 
B''\P(\lnorm Q(f,\bm\xi^{k})\rnorm\ge t). 
\]
\end{itemize}
\end{theorem}
\Proof

{\bm (A$''$)}: By symmetry, using Theorem \ref{contrRad}.(ii), with 
$\eta=\xi_1+\ldots+
\xi_k$ and $A=k$, we obtain that
\[
\P(Q(f, \bm\xi^k)\rnorm \ge t K) \le 
 K\P[ k^k Q( f, (\bm\xi_1+\ldots+\bm\xi_k)^k)\rnorm \ge t K)
\]
By \refrm{newinter}, and inequality \cite[(6.9.5)]{KwaW:book}, 
the latter quantity
can be estimated from below by
\[
c_k\,K\,
\P(k^{2k} k^k Q(f,  \bm\xi_1,\ldots,\bm\xi_k)\rnorm \ge t),
\]
which completes the proof of (A$''$).

{\bf (B$''$)}: By the 
Mazur-Orlicz polarization formula \refrm{MO}, and \refrm{inter}, 
we obtain the following estimates
\[
\begin{array}{rl}
\P(\lnorm
\wh{Q}(f,\bm\xi_1,\ldots,\bm\xi_k)\rnorm \ge t)
=&\P(\lnorm  \Rec{k!}Q( f,\Sum_{\ddl} 
(-1)^{k-|\ddl|} 
(\dl_1\bm\xi_1+\ldots+\dl_k\bm\xi_k)^k)\rnorm \ge t)\\
\le &\Sum_{\ddl}
\P(2^k\lnorm  {\Rec{k!}} Q(f, (\dl_1\bm\xi_1+\ldots+\dl_k\bm\xi_k)^k)
\rnorm\ge t) \\
=&\displaystyle 
\Sum_{i=0}^k{k \choose i}
\P(2^k\lnorm  \Rec{k!}Q(f, (\bm\xi_1+\ldots+\bm\xi_i)^k) \rnorm\ge t). \\
\end{array}
\]
By  \refrm{newinter}, with $j_1=\ldots=j_k$,  and the inequality 
\cite[(6.9.5)]{KwaW:book}, we estimate the above expression from 
above by  
\[
\begin{array}{rl}
 &\displaystyle
c_k^{-1} \Sum_{i=0}^k{k \choose i}
\lnorm\P(\Frac{(2i)^k}{k!} Q( f,\bm\xi_1^k)\rnorm\|\ge t) \\
\le  &
c_k^{-1}2^k \P(\Frac{(2k)^k}{k!}\lnorm Q(f,\bm\xi_1^k) \rnorm\ge t), \\
\end{array}
\]
which completes the proof.
\QED

\begin{remark}\label{ABtail:rem}\rm
While the symmetry
assumption is irrelevant in condition (B$''$) (or in (B), before), the
symmetrization procedure used in the proof of (A) fails. 
The reason is, that we use the conditioning on $\ca G$, which 
destroys the independence, which is essential in applications 
of \refrm{centdec}. 
\end{remark}

\section{Notes}
{~~~~~~~~}
\begin{enumerate}
\item
The inverse estimate in \refrm{centdec} is not true, in general, even if $k=1$.
For example, let $\xi_1,\ldots,\xi_n$ be Bernoulli random variables with
$p=\P(\xi_1=1)=1/2$, and $f(i)=1$. Then $\E|f\ov{\bm\xi}|^2= n/4$ and
$\E|f\bm\xi|^2= (n+n^2)/4$.
\item
The symmetry of functions $f$ is essential in 
Theorem \ref{AB}.(B) and its analogs, as was pointed
out in \cite{McCT:bandec}. The Bourgain's counterexample, given there, involves
$\BE=\ell^2\otimes\ell^2$ endowed with the projective norm $\|\bm
a\|=\inf\set{\sum_{i,j}\|a^1_i\|\cdot\|a^2_j\|:\bm a= \sum_{i,j} a^1_i\otimes
a^2_j}$, Rademacher chaoses, and tetrahedral functions $f$.
However, the inequalities (B) of both
Theorems \ref{AB} and \ref{F} hold for tetrahedral Rademacher
chaoses induced by $\bm\xi$ and $\BX$ (with  independent columns), 
whenever $\BE$ is (a) a Banach lattices with no subspace isomorphic 
to $c_0$, or (b) a UMD-space.
\item
The full analog of
Theorem \ref{AB} is valid in locally convex spaces.
\item
The decoupling results from Section~2 can be carried over to 
linear spaces over the
field of complex numbers.  To obtain similar results for Section~3 is
more difficult.  One approach is to show that if $\beps$\ denotes a sequence
of independent Rademacher random variables, and if $\bm \sigma$\ denotes a
sequence of independent Steinhaus random variables (that
is, $\sigma_i$\ is uniformly distributed over the complex unit circle), then
$
\| Q(f; \beps^k) \| \approx \| Q(f; \bm \sigma^k)\| 
$.
We omit the details of the development.
\item
In the case when the tail decoupling holds, i.e., in  Theorems \ref{ABtail}, 
\ref{F}.(A$'''$) and (B$'''$), we obtain the comparison of tightness. That 
is, for a family of functions $\set{f:f\in F}$, we have that, if one type 
of chaos $\set{Q_d(f):f\in F}$ is tight, so is the other, 
$\set{\wt{Q_d}(f):f\in F}$, subject to restrictions listed in the above 
theorems. That remark also applies to functions $f$ taking values in a 
locally convex space.
\item
In the context discussed above, we immediately obtain the comparison of 
generalized Orlicz modulars, i.e., functionals of the form 
$\Phi(\cdot)=\E\phi(\|\cdot\|)$, where $\phi$ is a nondecreasing 
function on the positive half-line, $\phi(0)=0$.
\item 
Multiple stochastic integrals of deterministic multivariate functions 
(cf., e.g., \cite{KalS}) can be seen as limits of multilinear random 
forms. Therefore if $\bm\xi,\bm\xi_1,\bm\xi_k$ are stochastic 
processes with independent increments, and the symbols 
$\ll f\, \bm\xi_1\otimes\cdots\otimes\bm\xi_d\rr$ and 
$\ll f\,\bm\xi^{\otimes k}\rr$ are understood as such integrals, 
then all decoupling inequalities carry over word-for-word.
\item\label{riorl}
Our decoupling inequalities involve a certain means of 
domination. Essentially, we show that the domination 
by means of $L^p$-norms yields the same for probability tails. 
The passing from one to another type of domination may be of an 
intrinsic interest. Recall the definition of $f^{**}$\ mentioned in 
Section~\ref{lorentzeg}. Let us note the following 
result, which can be applied in a wider context than ours.  

Suppose that 
$\xi$ and $\eta$\ are two given non-negative random variables, and define
quantities $c_1,\dots,c_5$\ below.
\begin{enumerate}
\item[{\rm (i)}]  Let $c_1$\ be the smallest constant such that 
for every Orlicz function, 
$\|\xi\|_{\phi}\le c_1 \|\eta\|_{\phi}$;
\item[{\rm (ii)}] Let $c_2$\ be the smallest constant such that
for all $t>0$, if $\phi_t(x)=(x-1)_+/t$, then
$\|\xi\|_{\phi_t}\le c_2 \|\eta\|_{\phi_t}$;
\item[{\rm (iii)}] Let $c_3$\ be the smallest constant such that
$\xi^{**} \le c_3 \eta^{**}$;
\item[{\rm (iv)}] Let $c_4$\ be the smallest constant such that
for every r.i. norm, $\|\xi\|\le c_4 \|\eta\|$;
\item[{\rm (v)}] Let $c_5$\ be the smallest constant such that
for every separable r.i. norm, $\|\xi\|\le c_5 \|\eta\|$.
\end{enumerate}
Then $c_1 = c_2 \le c_3 = c_4 = c_5 \le 2 c_1$.  Indeed, inequalities
$c_2 \le c_1 \le
c_4$\ and $c_3 \le c_5 \le c_4$\ are obvious.  That $c_1 \le c_2$\ follows
immediately from the formula
\[ \phi(x) = \int_0^\infty \phi_t(x) \, d(\phi'(t)) .\]
That $c_4 \le c_3$\ 
was proved in \cite[Proposition 2.a.8]{LinT}. 
That $c_3 \le 2c_2$\ follows from the formula 
$\|\xi\|_{\phi_t} \le \xi^{**}(t) \le 2\|\xi\|_{\phi_t}$.
To show the left hand side,
suppose that $\xi^{**}(t) \le 1$.
Then
\[ \int_0^t \xi^*(s) \, ds \le t .\]
Thus we have that $\xi^*(t) \le 1$, and hence
\[ 
\E\phi_t(\xi) = {1\over t}\int_0^1 (\xi^*(s)-1)_+ \, ds =
   {1\over t}\int_0^t (\xi^*(s)-1)_+ \, ds \le {1\over t}\int_0^t \xi^*(s) 
   \, ds \le 1 .
\]
To show the right hand side, suppose that $\|\xi\|_{\phi_t}\le 1$.  
Thus
\[ \int_0^{t_0} (\xi^*(s) - 1) \, ds \le t ,\]
where $t_0 = \P(\xi > 1)$.  If $t_0 \ge t$, then it follows that
\[ \int_0^{t} (\xi^*(s) - 1) \, ds \le t ,\]
from whence it follows that
\[ \int_0^{t} \xi^*(s) \, ds \le 2t .\]
If $t_0 < t$, then
\[ \int_0^t \xi^*(s) \, ds = \int_0^{t_0} (\xi^*(s) - 1) \, ds
   + \int_0^{t_0} \, ds + \int_{t_0}^t \xi^*(s) \, ds
   \le 2t ,\]
because $\xi^*(s) \le 1$\ if $s > t_0$.
\item
A decoupling principle for multivalued functions (proved in \cite{Pen})
also follows from our 
basic decoupling inequalities. Suppose that
$F(\cdot,\bm\xi)$ is a countably multivalued function, i.e., a countable family
of functions $\ca F_{\bi}$ is associated with each $\bi$.  In equivalent terms,
one may think of a decision function $\tau:D\times \N^k\to
\prod_{\bi\in\N^k}\ca F_{\bi}$ ($D$ is countable). Then the statements of
Theorem \ref{F} hold uniformly with respect to $\tau$, that is, the norm
$\|F(\cdot)\|$ is replaced by $\sup_\tau\sup_d\|\tau(d,\cdot)(\cdot)\|$. 
The theorem follows for a finite collection of decision functions
$\set{\tau_1,\ldots,\tau_n}$, since this means the replacement of
the underlying Banach space $\BE$ by another Banach space $\ell_n^\infty(\BE)$.
In the full statement we need the banach lattice $\L(\ell^\infty)$ to 
satisfy the property ``{\em $\sup_n\|x_n\|=\|\sup _n x_n\|$, for an 
increasing sequence of 
nonnegative vectors}''.  In view of the preceding note, we may choose  
a family of Orlicz spaces, and the required property holds.

Other sequential functionals on $\R^D$, e.g.
$\ell^p$, Orlicz $\ell^\psi$, etc., yield numerous variations of Theorem
\ref{F}.
\end{enumerate}

\begin{remark}\rm
 This paper represents the combination of the papers
\cite{PenM} and \cite{Szu:dechom}
\end{remark}

\begin{thebibliography}{dlPnMS92}

\bibitem[AMS92]{AsmM}
N.H. Asmar and S.J. Montgomery-Smith.
\newblock On the distribution of {Sidon} series.
\newblock {\em Arkiv f\"or Mat.}, 31:13--26, 1993.

\bibitem[Bur86]{Bur:mf}
D.L. Burkholder.
\newblock {\em Martingales and {Fourier} analysis in {Banach} space}, 
{\em \LNM}   1206,  61--108.
\newblock Springer Vg., C.I.M.E. Lectures, Varenna, Italy, 1985, 1986.

\bibitem[DA87]{DeA:dec}
A.~De~Acosta.
\newblock A decoupling inequality for multilinear forms of stable vectors.
\newblock {\em \PMS}, 8:71--76, 1987.

\bibitem[dlPn92]{Pen}
V.H. de~la Pe\~{n}a.
\newblock Decoupling and {Khintchine}'s inequalities for {U}-statistics.
\newblock {\em \AP}, 20:1877--1892, 1992.

\bibitem[dlPn93]{Pen:Wald}
V.H. de~la Pe\~{n}a.
\newblock Inequalities for tails of adapted process with an application to
  {W}ald's {L}emma.
\newblock {\em J.\ Theoretical Prob.}, 6:~, 1993.

\bibitem[dlPn94]{Pen:bound}
V.H. de~la Pe\~{n}a.
\newblock A bound on the moment generating function of a sum of dependent
  variables with an application to simple random sampling without replacement.
\newblock to appear in {\em Ann.\ Inst.\ Henri Poincar\'e}, 1994.

\bibitem[dlPnMS92]{PenM}
V.H. de~la Pe\~na and S.J. Montgomery-Smith.
\newblock Decoupling inequalities for tail probabilities of multilinear forms
  in symmetric and hypercontractive variables.
\newblock Preprint, 1992.

\bibitem[Doo53]{Doo}
Doob, J. (1953).
\newblock {\em Stochastic Processes}.
\newblock Wiley, New York-London-Sydney.

\bibitem[Hit88]{Hit:tan}
P.~Hitczenko.
\newblock Comparison of moments for tangent sequences of random variables.
\newblock {\em \PTRF}, 78:223--230, 1988.

\bibitem[KS88]{KraS:hyp}
W.~Krakowiak and J.~Szulga.
\newblock Hypercontraction principle and random multilinear forms in {Banach}
  spaces.
\newblock {\em \PTRF}, 77:325--342, 1988.

\bibitem[KS89]{KalS}
O.~Kallenberg and J.~Szulga.
\newblock Multiple integration with respect to {Poisson} and {L\'evy}
  processes.
\newblock {\em \PTRF}, 83:101--134, 1989.

\bibitem[KW87]{KwaW:double}
S.\ Kwapie\'n and W.A. Woyczy\'nski.
\newblock Double stochastic integrals, random quadratic forms and random series
  in {Orlicz} spaces.
\newblock {\em \AP}, 15:1072--1096, 1987.

\bibitem[KW92]{KwaW:book}
S.~Kwapie\'n and W.A. Woyczy\'nski.
\newblock {\em Random series and stochastic integrals}.
\newblock B{\"u}rkhauser, Boston, 1992.

\bibitem[Kwa87]{Kwa:dec}
S.~Kwapie\'n.
\newblock Decoupling inequalities and polynomial chaos.
\newblock {\em \AP}, 15:1062--1071, 1987.

\bibitem[LT79]{LinT}
J.\ Lindenstrauss and L.\ Tzafriri.
\newblock {\em Classical {Banach} spaces {II}. Function Spaces}.
\newblock \Spr, 1979.

\bibitem[LT91]{LedT}
M.~Ledoux and M.~Talagrand.
\newblock {\em Probability in Banach spaces}.
\newblock \Spr, Berlin, 1991.

\bibitem[MO35]{MazO}
S.~Mazur and W.~Orlicz.
\newblock Grundlegende {Eigenschaften} der polynomischen {Operationen}.
\newblock {\em \SM}, 5:50--68, 179--189, 1935.

\bibitem[MT86a]{McCT:realdec}
T.R. Mcconnell and M.S. Taqqu.
\newblock Decoupling inequalities for multilinear forms in independent
  symmetric random variables.
\newblock {\em \AP}, 14:943--954, 1986.

\bibitem[MT86b]{McCT:dyadic}
T.~McConnell and M.~Taqqu.
\newblock Dyadic approximation of double integrals with respect to symmetric
  stable processes.
\newblock {\em Stochastic Processes and Appl.}, 22:323--331, 1986.

\bibitem[MT87]{McCT:bandec}
T.R. Mcconnell and M.S. Taqqu.
\newblock Decoupling of {Banach}-valued multilinear forms in independent
  symmetric {Banach}-valued random variables.
\newblock {\em \PTRF}, 75:499--507, 1987.

\bibitem[MZ37]{MarZ}
J.~Marcinkiewicz and A.~Zygmund.
\newblock Sur les fonctions independantes.
\newblock {\em Fund.\ Math.}, 29:60--90, 1937.

\bibitem[NP87]{NolP}
D.~Nolan and D.~Pollard.
\newblock {U}-processes: rates and convergence.
\newblock {\em Annals of Stat.}, 15(2):780--799, 1987.

\bibitem[PZ32]{PalZ}
R.E.A.C. Paley and A.~Zygmund.
\newblock A note on analytic functions on the circle.
\newblock {\em Proc.\ Cambridge Phil. Soc.}, 28:266--272, 1932.

\bibitem[RST91]{RST}
J.\ Rosi\'nski, G.\ Samorodnitsky, and M.S. Taqqu.
\newblock Sample path properties of stochastic processes represented as
  multiple stable integrals.
\newblock {\em J. Multivariate Analysis}, 37:115--134, 1991.

\bibitem[RW86]{RosW:clock}
J.~Rosi\'nski and W.A. Woyczy\'nski.
\newblock On {It\^o} stochastic integration with respect to p-stable motion:
  Inner clock, integrability of sample paths, double and multiple integrals.
\newblock {\em \AP}, 14:271--286, 1986.

\bibitem[Szu91]{Szu:Umli}
J.~Szulga.
\newblock Limit theorems for some randomized nonlinear functionals of empirical
  measures.
\newblock \AU, Preprint, 1991.

\bibitem[Szu92]{Szu:dechom}
J.~Szulga.
\newblock Robust decoupling of homogeneous random chaoses.
\newblock \AU, Preprint, 1992.

\bibitem[{Wie}30]{Wie:chaos}
N.~{Wiener}.
\newblock The homogeneous chaos.
\newblock {\em \AJM}, 60:897--936, 1930.

\bibitem[WW43]{WieW}
W.\ Wiener and A.~Wintner.
\newblock The discrete chaos.
\newblock {\em \AJM}, 65:279--298, 1943.

\bibitem[Zin86]{Zin:comp}
J.~Zinn.
\newblock Comparison of martingale difference sequences.
\newblock In A.~Beck et~al., editor, {\em Probability on {Banach} spaces},
  \LNM, pages 453--457. \Spr~ 1153, 1986.

\end{thebibliography}

\vfill~\hfill 
\begin{tabular}{r}
\begin{tabular}{|l|}\hline
\rule{3in}{0pt}\\
\small Victor H. de la Pe\~na\\
\small Department of Statistics, \\
\small Columbia University,  \\
\small New York, NY 10027\\
\begin{tabular}{rl}
\small phone &\small (212) 854-5360\\
\small email &vp@stat.columbia.edu\\
\end{tabular}\\ \hline
\end{tabular}
\\
\\
\begin{tabular}{|l|}\hline
\rule{3in}{0pt}\\
\small Stephen J.~Montgomery-Smith\\
\small Department of Mathematics, \\
\small University of Missouri, \\
\small Columbia, MO 65211\\
\begin{tabular}{rl}
\small phone &\small (314) 882-7492\\
\small email &stephen@mont.cs.missouri.edu\\
\end{tabular}\\ \hline
\end{tabular}
\\
\\
\begin{tabular}{|l|}\hline
\rule{3in}{0pt}\\
\small Jerzy Szulga\\
\small Department of Mathematics,  \\
\small Auburn University,  \\
\small Auburn, AL 36849-3501\\
\begin{tabular}{rl}
\small phone &\small (205) 844-6569\\
\small email &szulgje@mail.auburn.edu\\
\end{tabular}\\ \hline
\end{tabular}\\
\end{tabular}
\end{document}